\documentclass[11pt,british]{amsart}
\usepackage{ifpdf}
\ifpdf\relax\else
\usepackage[T1]{fontenc}
\fi

\usepackage{babel,eucal,url,amssymb,stmaryrd,
enumerate,amscd,
}

\usepackage[pagebackref]{hyperref}

\theoremstyle{plain}
\newtheorem{thm}{Theorem}[section]
\newtheorem{lem}{Lemma}[section]
\newtheorem{pro}{Proposition}[section]

\theoremstyle{definition}
\newtheorem{defn}{Definition}[section]

\theoremstyle{remark}
\newtheorem{rem}{Remark}[section]
\newtheorem{example}{Example}[section]

\newcommand{\Gtwo}{\ifmmode{{\rm G}_2}\else{${\rm G}_2$}\fi}

\DeclareMathOperator{\tr}{tr}

\date{\today}

\begin{document}

\title{ALMOST PARACONTACT MANIFOLDS}

\author{Galia Nakova}
\address[Nakova]{University of Veliko Turnovo "St. Cyrilend St. Metodius"\\
Faculty of Education,\\
T. Tarnovski 2 str.,\\
5003 Veliko Tarnovo, Bulgaria} \email{gnakova@yahoo.com}

\author{Simeon Zamkovoy}
\address[Zamkovoy]{University of Sofia "St. Kl. Ohridski"\\
Faculty of Mathematics and Informatics,\\
Blvd. James Bourchier 5,\\
1164 Sofia, Bulgaria} \email{zamkovoy@fmi.uni-sofia.bg}

\begin{abstract}
{In this paper eleven basic classes of almost paracontact
manifolds are introduced and some examples are constructed.}

\end{abstract}

\maketitle \setcounter{tocdepth}{3} \tableofcontents

\section*{Introduction}

As is well-known, in \cite{GH} almost Hermitian manifolds are
classified with respect to the decomposition in subspaces
invariant under the action of the structural group $U(n)$.Thus we
have an adequate framework for several types of almost Hermitian
manifolds, previously defined by a number of authors in terms of
geometric properties which retain some portion of K\"{a}hler
geometry. The previous method was used in \cite{N} for Riemannian
almost product manifolds, and in \cite{GB} for almost complex
manifolds with Norden metric.

The geometry of almost contact manifolds is a natural extension in
the odd dimensional case of almost Hermitian geometry. Similarly,
the geometry of almost contact manifolds with $B$-metric can be
considered as a natural extension in the odd dimensional case of
geometry of almost complex Riemannian. A classification of almost
contact manifolds with $B$-metric with respect to the covariant
derivative of the fundamental tensor of type $(1,1)$ is made in
\cite{GMG}. The authors obtain eleven basic classes of almost
contact manifolds with $B$-metric and construct some examples.

A classification of almost paraHermitian manifolds is made in
\cite{B}.

The authors give examples of the primitive classes, which are
based on general almost paraHermitian structure on the tangent
bundles given in \cite{C}.

A classification of the almost paracontact Riemannian manifolds of
type $(n,n)$ with respect to the covariant derivative of type
$(1,1)$-tensor of the almost paracontact structure is made in
\cite{MS}. The authors consider almost paracontact Riemannian
manifolds of type $(n,n)$ with positive definite Riemannian metric
$g$, which is compatible with almost paracontact structure and it
satisfies the condition $g(\varphi \cdot,\varphi
\cdot)=g(\cdot,\cdot)-\eta (\cdot)\eta (\cdot)$.

The method used in the present paper is analogous of the one used
in \cite{GMG}. We give a classification of the almost paracontact
manifolds with respect to the covariant derivative of the
$(1,1)$-tensor of the almost paracontact structure. We consider
almost paracontact pseudo-Riemannian manifolds with indefinite
metric $g$, which it compatible with almost paracontact structure
and it satisfies the condition \eqref{con}. We obtain eleven basic
classes and construct some examples.

\section{Preliminaries}\label{sec-prelim}
A (2n+1)-dimensional smooth manifold $M^{(2n+1)}$ has an
\emph{almost paracontact structure} $(\varphi,\xi,\eta)$ if it
admits a tensor field $\varphi$ of type $(1,1)$, a vector field
$\xi$ and a 1-form $\eta$ satisfying the  following compatibility
conditions \cite{KW,Zam}:
\begin{eqnarray}
  \label{f82}
    & &
    \begin{array}{cl}
          (i)   & \varphi(\xi)=0,\quad \eta \circ \varphi=0,\quad
          \\[5pt]
          (ii)  & \eta (\xi)=1 \quad \varphi^2 = id - \eta \otimes \xi,
          \\[5pt]
          (iii) & \textrm{let $\mathbb D=Ker~\eta$ be the  horizontal distribution generated
                          by $\eta$, then}
          \\[3pt]
                & \textrm{the tensor field $\varphi$ induces an almost paracomplex structure on}
          \\[3pt]
                &  \textrm{each fibre on $\mathbb D$.}
    \end{array}
\end{eqnarray}

Recall that an almost paracomplex structure on an 2n-dimensional
manifold is a (1,1)-tensor $J$ such that $J^2=1$ and the
eigensubbundles $T^+,T^-$ corresponding to the eigenvalues $1,-1$
of $J$, respectively have equal dimension $n$. The Nijenhuis
tensor $N$ of $J$, given by
$N_{J}(X,Y)=[JX,JY]-J[JX,Y]-J[X,JY]+[X,Y],$ is the obstruction for
the integrability of the eigensubbundles $T^+,T^-$. If $N=0$ then
the almost paracomplex structure is called paracomplex or
integrable.

An immediate consequence of the definition of the almost
paracontact structure is that the endomorphism $\varphi$ has rank
$2n$, $\varphi \xi=0$ and $\eta \circ \varphi=0$ \cite{Zam}.

If a manifold $M^{(2n+1)}$ with $(\varphi,\xi,\eta)$-structure
admits a pseudo-Riemannian metric $g$ such that
\begin{equation}\label{con}
g(\varphi X,\varphi Y)=-g(X,Y)+\eta (X)\eta (Y),
\end{equation}
then we say that $M^{(2n+1)}$ has an almost paracontact metric
structure and $g$ is called \emph{compatible} metric. Any
compatible metric $g$ with a given almost paracontact structure is
necessarily of signature $(n+1,n)$ \cite{Zam}.

Setting $Y=\xi$, we have $\eta(X)=g(X,\xi).$

The fundamental 2-form
\begin{equation}\label{fund}
F(X,Y)=g(\varphi X,Y)
\end{equation}
is non-degenerate on the horizontal distribution $\mathbb D$ and
$\eta\wedge F^n\not=0$.

We have the following \cite{Zam}
\begin{defn}
If $g(X,\varphi Y)=d\eta(X,Y)$ (where
$d\eta(X,Y)=\frac12(X\eta(Y)-Y\eta(X)-\eta([X,Y])$ then $\eta$ is
a paracontact form and the almost paracontact metric manifold
$(M,\varphi,\eta,g)$ is said to be $\emph{paracontact metric
manifold}$.
\end{defn}

For a manifold $M^{(2n+1)}$ with an almost paracontact metric
structure $(\varphi,\xi,\eta,g)$ we can also construct a useful
local orthonormal basis. Let U be a coordinate neighborhood on M
and $X_1$ any unit vector field on U orthogonal to $\xi$. Then
$\varphi X_1$ is a vector field orthogonal to  both X and $\xi$,
and $|\varphi X_1|^2=-1$. Now choose a unit vector field $X_2$
orthogonal to $\xi$, $X_1$ and $\varphi X_1$. Then $\varphi X_2$
is also vector field orthogonal to $\xi$, $X_1$, $\varphi X_1$ and
$X_2$, and $|\varphi X_2|^2=-1$. Proceeding in this way we obtain
a local orthonormal basis $(X_i,\varphi X_i,\xi),i=1...n$ called a
\emph{$\varphi$-basis}.

Hence, an almost paracontact metric manifold
$(M^{2n+1},\varphi,\eta,\xi,g)$ is an odd dimensional manifold
with a structure group $\mathbb U(n,\mathbb R)\times Id$, where
$\mathbb U(n,\mathbb R)$ is the para-unitary group isomorphic to
$\mathbb {GL}(n,\mathbb R)$.

A paracontact structure for which $\xi$ is Killing vector field is
called a \emph{K-paracontact structure}.
\par
Let $\nabla$ be the Levi-Civita connection of the compatible
metric $g$. For all vectors $X,Y,Z \in
T_pM$, $p \in M$, we denote
\begin{equation}\label{4}
F(X,Y,Z)=g((\nabla_X\varphi)Y,Z).
\end{equation}
From \eqref{f82} and \eqref{con} the tensor $F$ has the following
properties:
\begin{equation}\label{5}
\begin{array}{lll}
F(X,Y,Z)=-F(X,Z,Y), \\
F(X,\varphi Y, \varphi
Z)=F(X,Y,Z)+\eta(Y)F(X,Z,\xi)-\eta(Z)F(X,Y,\xi),
\end{array}
\end{equation}
for all vectors $X,Y,Z \in T_pM$.
\par
The following 1-forms are associated with $F$:
\begin{equation}\label{6}
\theta(X)=g^{ij}F(e_i,e_j,X), \quad
\theta^*(X)=g^{ij}F(e_i,\varphi e_j,X), \quad
\omega(X)=F(\xi,\xi,X),
\end{equation}
where $X \in T_pM$, $\{e_i,\xi\}$, $(i=1,\ldots,2n)$ is a basis of
$T_pM$, and $(g^{ij})$ is the inverse matrix of $(g_{ij})$.

\section{The space of covariant derivatives of the structure $\varphi$}\label{sec-2}
\setcounter{equation}{0}

Let $V$ be a $(2n+1)$-dimensional vector space with almost
paracontact structure $(\varphi, \xi, \eta)$ and metric $g$
satisfying \eqref{con}. For an arbitrary $X \in V$ we have
$\varphi^2X=X-\eta(X)\xi \Longleftrightarrow
X=\varphi^2X+\eta(X)\xi$. Hence $V$ admits a decomposition into a
direct sum of vector subspaces
\[
V=\mathbb{D}\oplus \{\xi\},
\]
where $\mathbb{D}=Ker\eta$, $\{\xi\}=(Im\eta)\xi$. Then for an
arbitrary $X \in V$ it follows $X=hX+\eta(X)\xi$, where $X \in
\mathbb{D}, \eta(X)\xi \in \{\xi\}$. Denoting the restrictions of
$g$ and $\varphi$ on $\mathbb{D}$ with the same letters we obtain
an $2n$-dimension almost paracomplex manifold
$(\mathbb{D},\varphi,g)$.
\par
Let $\{e_1,\ldots,e_{2n}\}$ be any basis of $\mathbb{D}$. Then
$\{e_1,\ldots,e_{2n},\xi\}$ is a basis of $V$ and for an arbitrary
$X\in V$ we have $X=X^ie_i+\eta(X)\xi, \, i=1,\ldots,2n$.
\par
We define the linear operators
\begin{equation}\label{2.1}
\begin{array}{ll}
{\mathcal{A}} _{e_i} : V \longrightarrow V: \quad X\longrightarrow
{\mathcal{A}} _{e_i}X, \quad i=1,\ldots,2n; \\ \\
{\mathcal{A}} _\xi : V \longrightarrow \mathbb{D}: \quad X
\longrightarrow {\mathcal{A}} _\xi X,
\end{array}
\end{equation}
having the following properties
\begin{equation}\label{2.2}
g({\mathcal{A}} _{e_i}X,e_j)=-g({\mathcal{A}} _{e_j}X,e_i), \quad i,j=1,\ldots,2n;
\end{equation}
\begin{equation}\label{2.3}
{\mathcal{A}} _{\varphi e_i}X=-\varphi({\mathcal{A}} _{e_i}X)-g({\mathcal{A}} _\xi X,e_i)\xi;
\end{equation}
\begin{equation}\label{2.4}
\eta({\mathcal{A}} _{e_i}X)=-g({\mathcal{A}} _\xi X,\varphi e_i);
\end{equation}
\begin{equation}\label{2.5}
\eta({\mathcal{A}} _\xi X)=0.
\end{equation}
We consider the vector space ${\mathcal{F}}$ of all tensors $F$ of
type $(0,3)$ over $V$, defined by
\begin{equation}\label{2.6}
F(X,Y,Z)=Y^ig({\mathcal{A}} _{e_i}X,Z)+\eta(Y)g({\mathcal{A}} _\xi
X,\varphi Z),
\end{equation}
where ${\mathcal{A}} _{e_i} \, (i=1,\ldots,2n)$ and ${\mathcal{A}}
_\xi$ have the properties \eqref{2.2}$\div$\eqref{2.5}. It is easy
to verify that the tensors $F$ do not depend on the basis of $V$.
Using \eqref{2.6} and \eqref{2.2}$\div$\eqref{2.5} we establish
that the tensors $F \in {\mathcal{F}}$ have the properties
\eqref{5}.
\par
The compatible metric $g$ induces on ${\mathcal{F}}$ an inner
product $<,>$, defined by
\begin{equation}\label{2.7}
<F_1,F_2>=g^{iq}g^{jr}g^{ks}F_1(f_i,f_j,f_k)F_2(f_q,f_r,f_s)
\end{equation}
for $F_1, F_2 \in  {\mathcal{F}}$ and $\{f_1,\ldots,f_{2n+1}\}$ is
a basis of $V$.
\par
The standard representation of the structure group $\mathbb
U(n,\mathbb R)\times Id$ in $V$ induces a natural representation
$\lambda $ of $\mathbb U(n,\mathbb R)\times Id$ in
${\mathcal{F}}$:
\[
(\lambda (a)F)(X,Y,Z)=F(a^{-1}X,a^{-1}Y,a^{-1}Z), \quad a\in
\mathbb U(n,\mathbb R)\times Id, \quad F\in {\mathcal{F}}, \, \, X,Y,Z \in V,
\]
so that
\[
<\lambda (a)F_1,\lambda (a)F_2>=<F_1,F_2>, \quad a\in
\mathbb U(n,\mathbb R)\times Id, \quad F_1, F_2 \in {\mathcal{F}}.
\]
\par
Let $X \in V$. Then for ${\mathcal{A}} _{e_i}X \in V, \,
(i=1,\ldots,2n)$ we have
\begin{equation}\label{2.8}
{\mathcal{A}} _{e_i}X=h({\mathcal{A}}
_{e_i}hX)+\eta(X)h({\mathcal{A}} _{e_i}\xi)+\eta({\mathcal{A}}
_{e_i}hX)\xi+ \eta(X)\eta({\mathcal{A}} _{e_i}\xi)\xi .
\end{equation}
Taking into account \eqref{2.6} and \eqref{2.8} we obtain
\begin{equation}\label{2.9}
\begin{array}{ll}
F(X,Y,Z)=Y^ig(h({\mathcal{A}}
_{e_i}hX),Z)-Z^i\eta(Y)\eta({\mathcal{A}} _{e_i}hX)+Y^i\eta(Z)
\eta({\mathcal{A}} _{e_i}hX) \\ \\
+\eta(X)Y^ig(h({\mathcal{A}}
_{e_i}\xi),Z)+\eta(X)\eta(Z)Y^i\eta({\mathcal{A}} _{e_i}\xi)-
\eta(X)\eta(Y)Z^i\eta({\mathcal{A}} _{e_i}\xi).
\end{array}
\end{equation}
Analogously as \cite{GMG} we define the operators
\[
\begin{array}{l}
p_i: {\mathcal{F}} \longrightarrow {\mathcal{F}}; \quad i=1,2,3,4, \\ \\
p_1(F)(X,Y,Z)=F(hX,hY,hZ); \\ \\
p_2(F)(X,Y,Z)=-\eta(Y)F(hX,hZ,\xi)+\eta(Z)F(hX,hY,\xi); \\ \\
p_3(F)(X,Y,Z)=\eta(X)F(\xi,hY,hZ); \\ \\
p_4(F)(X,Y,Z)=\eta(X)\eta(Y)F(\xi,\xi,hZ)-\eta(X)\eta(Z)F(\xi,\xi,hY).
\end{array}
\]
\begin{lem}\label{l 2.1}
The operators $p_i \, (i=1,2,3,4)$ have the following properties
\[
\begin{array}{cl}
          (i)   & p_i\circ p_i=p_i, \quad i=1,2,3,4; \\ \\

          (ii)  & \displaystyle \sum^{4}_{i=1} p_i=id; \\ \\

          (iii) & p_i\circ p_j=0, \quad i \neq j, \quad i,j=1,2,3,4;
\end{array}
\]
and $p_i \, (i=1,2,3,4)$ commute with $\mathbb U(n,\mathbb R)\times Id$.
\end{lem}
We denote $W_i=Im p_i \quad (i=1,2,3,4)$.
\begin{pro}\label{p 2.1} ( Partial decomposition ) The decomposition
\[
{\mathcal{F}}=W_1\oplus W_2\oplus W_3\oplus W_4
\]
is orthogonal and invariant under the action of the
group $\mathbb U(n,\mathbb R)\times Id$.
\end{pro}
\par
{\bf Proof.} From well known algebraic result and Lemma~\ref{2.1} we obtain
the decomposition ${\mathcal{F}}=W_1\oplus W_2\oplus W_3\oplus W_4$ that is
$\mathbb U(n,\mathbb R)\times Id$ - invariant. By direct
computations using \eqref{2.7} we check that $W_i \bot W_j, \quad
i \neq j, \quad i,j=1,2,3,4$.
\par
From \eqref{2.9}, \eqref{2.4} by explicit calculations we have
\begin{equation}\label{2.10}
\begin{array}{ll}
p_1(F)(X,Y,Z)=F(hX,hY,hZ)=Y^ig(h({\mathcal{A}}_{e_i}hX),Z); \\  \\
p_2(F)(X,Y,Z)=-\eta(Y)F(hX,hZ,\xi)+\eta(Z)F(hX,hY,\xi)=\\
-Z^i\eta(Y)\eta({\mathcal{A}} _{e_i}hX)+Y^i\eta(Z)\eta({\mathcal{A}} _{e_i}hX)= \\ \\
\eta(Y)g({\mathcal{A}} _\xi hX,\varphi Z)-\eta(Z)g({\mathcal{A}} _\xi hX,\varphi Y);\\ \\
p_3(F)(X,Y,Z)=\eta(X)F(\xi,hY,hZ)=\eta(X)Y^ig(h({\mathcal{A}} _{e_i} \xi),Z); \\ \\
p_4(F)(X,Y,Z)=\eta(X)\eta(Y)F(\xi,\xi,hZ)-\eta(X)\eta(Z)F(\xi,\xi,hY)=\\ \\
-\eta(X)\eta(Y)Z^i\eta({\mathcal{A}} _{e_i}\xi)+
\eta(X)\eta(Z)Y^i\eta({\mathcal{A}} _{e_i}\xi),
\end{array}
\end{equation}
where $X,Y,Z \in V$ and $(i=1,\ldots,2n)$. Using \eqref{2.8},
\eqref{2.9}, \eqref{2.10} we obtain
\begin{pro}\label{p 2.2}
Let ${\mathcal{A}} _{e_i} \, (i,\ldots,2n)$ be the linear
operators, defined by \eqref{2.1} and having the properties
\eqref{2.2}, \eqref{2.3}, \eqref{2.4}, \eqref{2.5}. Then for an
arbitrary $F \in {\mathcal{F}}$ and $X,Y,Z \in V$ we have
\[
\begin{array}{clll}
          (i)   & {\mathcal{A}} _{e_i}X=h({\mathcal{A}} _{e_i}hX), \, (i=1,\ldots,2n)\Longleftrightarrow F=p_1F; \\ \\
          (ii)  & {\mathcal{A}} _{e_i}X=\eta({\mathcal{A}} _{e_i}hX)\xi =
          -g({\mathcal{A}} _\xi hX,\varphi e_i)\xi , \,
          (i=1,\ldots,2n)\Longleftrightarrow F=p_2F; \\ \\
          (iii) & {\mathcal{A}} _{e_i}X=\eta(X)h({\mathcal{A}} _{e_i}\xi ), \,
          (i=1,\ldots,2n)\Longleftrightarrow F=p_3F; \\ \\
          (iiii) & {\mathcal{A}} _{e_i}X=\eta(X)\eta({\mathcal{A}} _{e_i}\xi )\xi =
          -\eta(X)g({\mathcal{A}} _\xi \xi,\varphi e_i)\xi, \,
          (i=1,\ldots,2n)\Longleftrightarrow F=p_4F .
\end{array}
\]
\end{pro}
\section{The subspace $W_1$}\label{sec-3}
\setcounter{equation}{0}
From Proposition~\ref{p 2.2} we have
\[
W_1=\left\{F\in {\mathcal{F}}: F=p_1F\Longleftrightarrow
{\mathcal{A}} _{e_i}X=h({\mathcal{A}} _{e_i}hX) \, \, (i=1,\ldots,2n)\right\}.
\]
The condition ${\mathcal{A}} _{e_i}X=h({\mathcal{A}} _{e_i}hX) \, \,
(i=1,\ldots,2n)$ is
equivalent to
\begin{equation}\label{3.1}
\begin{array}{l}
{\mathcal{A}} _{e_i}\xi =0; \quad  \eta({\mathcal{A}} _{e_i}X)=0 \, \,(i=1,\ldots,2n).
\end{array}
\end{equation}
From equalities \eqref{2.4} and \eqref{3.1} we obtain
\begin{equation}\label{3.2}
\begin{array}{l}
{\mathcal{A}} _{\xi }X=0.
\end{array}
\end{equation}
Then the decomposition of $W_1$ over $V$ coincides with the decomposition
of $W$ over $D$, where the vector space $W$ is defined by
\[
W=\left\{F\in {\mathcal{F}}: F(X,Y,Z)=-F(X,Z,Y)=F(X,\varphi Y, \varphi Z),\, \,
X,Y,Z\in D\right\}.
\]
Taking into account \eqref{2.2}$\div$\eqref{2.5}, \eqref{2.9}, \eqref{3.1},
\eqref{3.2} we have
\begin{equation}\label{3.3}
\begin{array}{l}
W=\left\{F\in {\mathcal{F}}: F(X,Y,Z)=Y^ig({\mathcal{A}} _{e_i}X,Z):
g({\mathcal{A}} _{e_i}X,e_j)=-g({\mathcal{A}} _{e_j}X,e_i),
\right.
\end{array}
\end{equation}
\[
\left.
{\mathcal{A}} _{\varphi e_i}X=-\varphi ({\mathcal{A}} _{e_i}X), \, \,
(i=1,\ldots,2n); \, \, X,Y,Z\in D\right\}.
\]
Using \eqref{3.3} and \eqref{6} we find
\[
\theta(Z)=-Z^i{\tr}{\mathcal{A}} _{e_i}, \quad \theta ^*(Z)=\theta(\varphi Z)=
Z^i{\tr}({\mathcal{A}} _{e_i}\circ\varphi ), \quad (i=1,\ldots,2n).
\]
We define the operators
\[
\begin{array}{ll}
m_i: W \longrightarrow W, \quad i=1,2; \\
m_1(F)(X,Y,Z)=\displaystyle{\frac{1}{2}}Y^i\left\{g({\mathcal{A}} _{e_i}X,
Z)-g({\mathcal{A}} _{e_i}\varphi X,\varphi Z)\right\} \, \,
(i=1,\ldots,2n); \\ \\
m_2(F)(X,Y,Z)=\displaystyle{\frac{1}{2}}Y^i\left\{g({\mathcal{A}} _{e_i}X,
Z)+g({\mathcal{A}} _{e_i}\varphi X,\varphi Z)\right\} \, \, (i=1,\ldots,2n).
\end{array}
\]
\begin{lem}\label{l 3.1}
The operators $m_i \, (i=1,2)$ have the following properties
\[
\begin{array}{cl}
          (i)   & m_i\circ m_i=m_i, \quad i=1,2; \\ \\

          (ii)  & \displaystyle \sum^{2}_{i=1} m_i=id; \\ \\

          (iii) & m_1\circ m_2=m_2\circ m_1=0;
\end{array}
\]
and $m_i \, (i=1,2)$ commute with $\mathbb U(n,\mathbb R)$.
\end{lem}
We denote $W_{11}=Im \, m_1$, ${\mathcal{F}}_3=Im \, m_2$. Lemma~\ref{l 3.1} implies
the decomposition
$W=W_{11}\oplus {\mathcal{F}}_3$ that is $\mathbb U(n,\mathbb R)$ - invariant.
Using \eqref{2.7} we check that $W_{11} \bot {\mathcal{F}}_3$.
\begin{pro}\label{p 3.1} The decomposition
\[
W=W_{11}\oplus {\mathcal{F}}_3
\]
is orthogonal and invariant under the action of the
group $\mathbb U(n,\mathbb R)$.
\end{pro}
After direct computations using definitions of $m_1, m_2$,
\eqref{3.3} we obtain
\begin{pro}\label{p 3.2}
For an arbitrary $F \in W$ we have
\[
\begin{array}{cl}
          (i)   & F=m_1F \Longleftrightarrow {\mathcal{A}} _{e_i} \circ \varphi =
          \varphi \circ {\mathcal{A}} _{e_i}, \, \, (i=1,\ldots,2n); \\
          (ii)  & F=m_1F \Longleftrightarrow F(X,Y,Z)=-F(\varphi X,\varphi Y,Z); \\
          (iii) & F=m_2F \Longleftrightarrow {\mathcal{A}} _{e_i} \circ \varphi =-
          \varphi \circ {\mathcal{A}} _{e_i}, \, \, (i=1,\ldots,2n); \\
          (iiii)& F=m_2F \Longleftrightarrow F(X,Y,Z)=F(\varphi X,\varphi Y,Z).
\end{array}
\]
\end{pro}
From Proposition~\ref{p 3.2} the characteristic conditions of
$W_{11}, {\mathcal{F}}_3$ are
\begin{equation}\label{3.4}
\begin{array}{ll}
W_{11}=\left\{F \in W: {\mathcal{A}} _{e_i} \circ \varphi =
\varphi \circ {\mathcal{A}} _{e_i}, \, \, (i=1,\ldots,2n)\right\}\Longleftrightarrow \\ \\
\left\{F \in W: F(X,Y,Z)=-F(\varphi X,\varphi Y,Z)\right\};
\end{array}
\end{equation}
\begin{equation}\label{3.5}
\begin{array}{ll}
{\mathcal{F}}_3=\left\{F \in W: {\mathcal{A}} _{e_i} \circ \varphi =-
\varphi \circ {\mathcal{A}} _{e_i}, \, \, (i=1,\ldots,2n)\right\}\Longleftrightarrow \\ \\
\left\{F \in W: F(X,Y,Z)=F(\varphi X,\varphi Y,Z)\right\}.
\end{array}
\end{equation}
We define the operator
\[
\begin{array}{ll}
m_3: W_{11} \longrightarrow W_{11}; \\
m_3(F)(X,Y,Z)=F(X,Y,Z)-\displaystyle{\frac{1}{2(n-1)}}\left\{g(X,\varphi Y)
\theta _{F}(\varphi Z)-g(X,\varphi Z)\theta _{F}(\varphi Y)-\right.
\end{array}
\]
\[
\left.g(\varphi X,\varphi Y)
\theta _{F}(Z)+g(\varphi X,\varphi Z)\theta _{F}(Y)\right\}.
\]
\begin{lem}\label{l 3.3}
The operator $m_3$ has the following properties
\[
\begin{array}{cl}
          (i)   & m_3\circ m_3=m_3; \\ \\

          (ii)  &  <m_3F_1,F_2>=<F_1,m_3F_2>, \quad F_1, F_2 \in W_{11}; \\ \\
\end{array}
\]
and $m_3$ commutes with $\mathbb U(n,\mathbb R)$.
\end{lem}
If we denote ${\mathcal{F}}_1=Ker m_3$ and ${\mathcal{F}}_2=Im \, m_3$,
then Lemma~\ref{l 3.3} implies
\begin{pro}\label{p 3.5} The decomposition
\[
W_{11}={\mathcal{F}}_1 \oplus {\mathcal{F}}_2
\]
is orthogonal and invariant under the action of the
group $\mathbb U(n,\mathbb R)$, \\ where
\[
\begin{array}{l}
{\mathcal{F}}_1=\{F \in W:
F(X,Y,Z)=\displaystyle{\frac{1}{2(n-1)}}\left(g(X,\varphi
Y)\theta(\varphi Z)-g(X,\varphi Z)\theta(\varphi Y)\right.
\end{array}
\]
\[
\left.\left.-g(\varphi X,\varphi Y)\theta(Z)+g(\varphi X,\varphi Z)\theta(Y)
\right)\right\},
\Longleftrightarrow
\]
\[
\begin{array}{l}
{\mathcal{F}}_1=\left\{F \in W: {\mathcal{A}}_{e_i}=
\displaystyle{\frac{{\tr}{\mathcal{A}} _{e_i}}{2n}} \, id+
\displaystyle{\frac{{\tr}({\mathcal{A}} _{e_i}\circ\varphi )}{2n}}
\, \varphi , \quad (i=1,\ldots,2n)\right\},
\end{array}
\]
\[
\begin{array}{l}
{\mathcal{F}}_2=\left\{F \in W: F(X,Y,Z)=-F(\varphi X,\varphi Y,Z), \quad
\theta =0\right\}\Longleftrightarrow
\end{array}
\]
\[
\begin{array}{l}
{\mathcal{F}}_2=\left\{F \in W: {\mathcal{A}} _{e_i}\circ\varphi
=\varphi\circ {\mathcal{A}} _{e_i},
{\tr}({\mathcal{A}}
_{e_i}\circ\varphi )=0, \quad (i=1,\ldots,2n)\right\}.
\end{array}
\]
\end{pro}
Taking into account Proposition~\ref{p 3.1}, Proposition~\ref{p
3.5} we obtain
\begin{pro}\label{3.4}
The decomposition
\[
W_1={\mathcal{F}}_1 \oplus {\mathcal{F}}_2 \oplus {\mathcal{F}}_3
\]
is orthogonal and invariant under the action of the
group $\mathbb U(n,\mathbb R)\times Id$.
\end{pro}

\section{The subspace $W_2$}\label{sec-4}
\setcounter{equation}{0} From $W_2=\{F \in {\mathcal{F}}:
F=p_2F\}$, \eqref{2.10}, Proposition~\ref{p 2.2} and \eqref{2.4}
it follows
\begin{equation}\label{4.111}
\begin{array}{ll}
W_2=\left\{F \in {\mathcal{F}}: F(X,Y,Z)=-\eta(Y)g(\varphi ({\mathcal{A}} _\xi X),Z)+ \right. \\ \\
\left.
\eta(Z)g(\varphi ({\mathcal{A}} _\xi X),Y), \, {\mathcal{A}} _\xi \xi =0\right\}.
\end{array}
\end{equation}
Using \eqref{4.111} and \eqref{6} we find
\[
\theta(\xi)={\tr}({\mathcal{A}}_\xi \circ\varphi ), \quad \theta^*(\xi )=
-{\tr}{\mathcal{A}}_\xi .
\]
We define the operators
\[
\begin{array}{llllll}
q_i: W_2 \longrightarrow W_2, \quad i=1,2; \\
q_1(F)(X,Y,Z)=-\displaystyle{\frac{1}{2}}\eta(Y)\left\{g(\varphi ({\mathcal{A}} _\xi X),Z)+g({\mathcal{A}} _\xi
(\varphi X),Z)\right\}+ \\ \\
\displaystyle{\frac{1}{2}}\eta(Z)\left\{g(\varphi ({\mathcal{A}} _\xi X),Y)+g({\mathcal{A}} _\xi
(\varphi X),Y)\right\}; \\ \\
q_2(F)(X,Y,Z)=-\displaystyle{\frac{1}{2}}\eta(Y)\left\{g(\varphi ({\mathcal{A}} _\xi X),Z)-g({\mathcal{A}} _\xi
(\varphi X),Z)\right\}+ \\ \\
\displaystyle{\frac{1}{2}}\eta(Z)\left\{g(\varphi ({\mathcal{A}} _\xi X),Y)-g({\mathcal{A}} _\xi
(\varphi X),Y)\right\}.
\end{array}
\]
\begin{lem}\label{l 4.1}
The operators $q_i \, (i=1,2)$ have the following properties
\[
\begin{array}{cl}
          (i)   & q_i\circ q_i=q_i, \quad i=1,2; \\ \\

          (ii)  & \displaystyle \sum^{2}_{i=1} q_i=id; \\ \\

          (iii) & q_1\circ q_2=q_2\circ q_1=0;
\end{array}
\]
and $q_i \, (i=1,2)$ commute with $\mathbb U(n,\mathbb R)\times Id$.
\end{lem}
We denote $W^{'}=Imq_1$, $W^{''}=Imq_2$. Lemma~\ref{l 4.1} implies
the decomposition
$W_2=W^{'}\oplus W^{''}$ that is $\mathbb U(n,\mathbb R)\times Id$ - invariant.
Using \eqref{2.7} we check that $W^{'} \bot W^{''}$.
\begin{pro}\label{p 4.1} The decomposition
\[
W_2=W^{'}\oplus W^{''}
\]
is orthogonal and invariant under the action of the
group $\mathbb U(n,\mathbb R)\times Id$.
\end{pro}
After direct computations using definitions of $q_1, q_2$,
\eqref{4.111}, \eqref{2.5} we obtain
\begin{pro}\label{p 4.2}
For an arbitrary $F \in W_2$ we have
\[
\begin{array}{cl}
          (i)   & F=q_1F \Longleftrightarrow {\mathcal{A}} _\xi \circ \varphi =
          \varphi \circ {\mathcal{A}} _\xi ; \\
          (ii)  & F=q_1F \Longleftrightarrow F(X,Y,Z)=-F(\varphi X,\varphi Y,Z)-F(\varphi X,Y,\varphi Z); \\
          (iii) & F=q_2F \Longleftrightarrow {\mathcal{A}} _\xi \circ \varphi =-
          \varphi \circ {\mathcal{A}} _\xi ; \\
          (iiii)& F=q_2F \Longleftrightarrow F(X,Y,Z)=F(\varphi X,\varphi Y,Z)+F(\varphi X,Y,\varphi Z).
\end{array}
\]
\end{pro}
From Proposition~\ref{p 4.2} the characteristic conditions of
$W^{'}, W^{''}$ are
\begin{equation}\label{4.2}
\begin{array}{ll}
W^{'}=\left\{F \in W_2: {\mathcal{A}} _\xi \circ \varphi =
\varphi \circ {\mathcal{A}} _\xi \right\}\Longleftrightarrow \\ \\
\left\{F \in W_2: F(X,Y,Z)=-F(\varphi X,\varphi Y,Z)-F(\varphi X,Y,\varphi Z)\right\};
\end{array}
\end{equation}
\begin{equation}\label{4.3}
\begin{array}{ll}
W^{''}=\left\{F \in W_2: {\mathcal{A}} _\xi \circ \varphi =-
\varphi \circ {\mathcal{A}} _\xi \right\}\Longleftrightarrow \\ \\
\left\{F \in W_2: F(X,Y,Z)=F(\varphi X,\varphi Y,Z)+F(\varphi X,Y,\varphi Z)\right\}.
\end{array}
\end{equation}
We define the operators
\[
\begin{array}{llllll}
r_i: W^{'} \longrightarrow W^{'}, \quad i=1,2; \\
r_1(F)(X,Y,Z)=-\displaystyle{\frac{1}{2}}\eta(Y)\left\{g(\varphi ({\mathcal{A}} _\xi X),Z)+g(\varphi X,
{\mathcal{A}} _\xi Z)\right\}+ \\ \\
\displaystyle{\frac{1}{2}}\eta(Z)\left\{g(\varphi ({\mathcal{A}} _\xi X),Y)+g(\varphi X,
{\mathcal{A}} _\xi Y)\right\}; \\ \\
r_2(F)(X,Y,Z)=-\displaystyle{\frac{1}{2}}\eta(Y)\left\{g(\varphi ({\mathcal{A}} _\xi X),Z)-g(\varphi X,
{\mathcal{A}} _\xi Z)\right\}+ \\ \\
\displaystyle{\frac{1}{2}}\eta(Z)\left\{g(\varphi ({\mathcal{A}} _\xi X),Y)-g(\varphi X,
{\mathcal{A}} _\xi Y)\right\}.
\end{array}
\]
\begin{lem}\label{l 4.2}
The operators $r_i \, (i=1,2)$ have the following properties
\[
\begin{array}{cl}
          (i)   & r_i\circ r_i=r_i, \quad i=1,2; \\ \\

          (ii)  & \displaystyle \sum^{2}_{i=1} r_i=id; \\ \\

          (iii) & r_1\circ r_2=r_2\circ r_1=0;
\end{array}
\]
and $r_i \, (i=1,2)$ commute with $\mathbb U(n,\mathbb R)\times Id$.
\end{lem}
We denote $W^{'}_1=Imr_1$, $W^{'}_2=Imr_2$. Lemma~\ref{l 4.2}
implies the decomposition $W^{'}=W^{'}_1\oplus W^{'}_2$ that is
$\mathbb U(n,\mathbb R)\times Id$ - invariant.
Using \eqref{2.7} we check that $W^{'}_1 \bot W^{'}_2$.
\begin{pro}\label{p 4.3} The decomposition
\[
W^{'}=W^{'}_1\oplus W^{'}_2
\]
is orthogonal and invariant under the action of the group
$\mathbb U(n,\mathbb R)\times Id$.
\end{pro}
Having in mind definitions of $r_1, r_2$, \eqref{4.1},
\eqref{4.2}, \eqref{2.5} we obtain
\begin{pro}\label{p 4.4}
For an arbitrary $F \in W^{'}$ we have
\[
\begin{array}{cl}
          (i)   & F=r_1F \Longleftrightarrow g({\mathcal{A}} _\xi \, . \, , \, . )=
          g( . \, , \, {\mathcal{A}} _\xi \, . ); \\
          (ii)  & F=r_1F \Longleftrightarrow F(X,Y,Z)=-F(Y,Z,X)+F(Z,X,Y)
          -2F(\varphi X,\varphi Y,Z); \\
          (iii) & F=r_2F \Longleftrightarrow
          g({\mathcal{A}} _\xi \, . \, , \, . )=-
          g( . \, , \, {\mathcal{A}} _\xi \, . ); \\
          (iiii)& F=r_2F \Longleftrightarrow F(X,Y,Z)=-F(Y,Z,X)-F(Z,X,Y).
\end{array}
\]
\end{pro}
From Proposition~\ref{p 4.4} the characteristic conditions of
$W^{'}_1, W^{'}_2$ are
\begin{equation}\label{4.4}
\begin{array}{ll}
W^{'}_1=\left\{F \in W^{'}: g({\mathcal{A}} _\xi \, . \, , \, . )=
g( . \, , \, {\mathcal{A}} _\xi \, . )\right\} \Longleftrightarrow \\ \\
\left\{F \in W^{'}: F(X,Y,Z)=-F(Y,Z,X)+F(Z,X,Y)-2F(\varphi X,\varphi Y,Z)\right\};
\end{array}
\end{equation}
\begin{equation}\label{4.5}
\begin{array}{ll}
W^{'}_2=\left\{F \in W^{'}: g({\mathcal{A}} _\xi \, . \, , \, . )=-
g( . \, , \, {\mathcal{A}} _\xi \, . )\right\} \Longleftrightarrow \\ \\
\left\{F \in W^{'}: F(X,Y,Z)=-F(Y,Z,X)-F(Z,X,Y)\right\}.
\end{array}
\end{equation}
We define the operator
\[
\begin{array}{lll}
s: W^{'}_1 \longrightarrow W^{'}_1; \\
s(F)(X,Y,Z)=F(X,Y,Z)+\displaystyle{\frac{\theta ^{*}_{F}(\xi )}{2n}}\left\{\eta(Y)g(X,\varphi Z)-
\eta(Z)g(X,\varphi Y)\right\}.
\end{array}
\]
\begin{lem}\label{l 4.3}
The operator $s$ has the following properties
\[
\begin{array}{cl}
          (i)   & s\circ s=s; \\ \\

          (ii)  &  <sF_1,F_2>=<F_1,sF_2>, \quad F_1, F_2 \in W^{'}_1; \\ \\
\end{array}
\]
and $s$ commutes with $\mathbb U(n,\mathbb R)\times Id$.
\end{lem}
If we denote ${\mathcal{F}}_5=Kers$ and ${\mathcal{F}}_6=Ims$,
then Lemma~\ref{l 4.3} implies
\begin{pro}\label{p 4.5} The decomposition
\[
W^{'}_1={\mathcal{F}}_5 \oplus {\mathcal{F}}_6
\]
is orthogonal and invariant under the action of the
group $\mathbb U(n,\mathbb R)\times Id$, \\ where
\[
\begin{array}{llll}
{\mathcal{F}}_5=\left\{F \in W^{'}_1: F(X,Y,Z)=-\displaystyle{\frac{\theta^*(\xi)}{2n}}\{\eta(Y)g(X,\varphi Z)-
\eta(Z)g(X,\varphi Y) \}\right\}\Longleftrightarrow \\ \\
\left\{F \in W^{'}_1: {\mathcal{A}}_\xi=\displaystyle{\frac{{\tr}
{\mathcal{A}}_\xi}{2n}} \, id\right\}, \\ \\
{\mathcal{F}}_6=\left\{F \in W^{'}_1: \theta^*(\xi)=0\right\}\Longleftrightarrow
\left\{F \in W^{'}_1: {\tr}{\mathcal{A}}_\xi=0\right\}.
\end{array}
\]
\end{pro}
We define the operator
\[
\begin{array}{lll}
t: W^{'}_2 \longrightarrow W^{'}_2; \\
t(F)(X,Y,Z)=F(X,Y,Z)-\displaystyle{\frac{\theta _{F}(\xi
)}{2n}}\left\{\eta(Y)g(\varphi X,\varphi Z)- \eta(Z)g(\varphi
X,\varphi Y)\right\}.
\end{array}
\]
\begin{lem}\label{l 4.4}
The operator $t$ has the following properties
\[
\begin{array}{cl}
          (i)   & t\circ t=t; \\ \\

          (ii)  &  <tF_1,F_2>=<F_1,tF_2>, \quad F_1, F_2 \in W^{'}_2; \\ \\
\end{array}
\]
and $t$ commutes with $\mathbb U(n,\mathbb R)\times Id$.
\end{lem}
If we denote ${\mathcal{F}}_4=Kert$ and ${\mathcal{F}}_7=Imt$,
then Lemma~\ref{l 4.4} implies
\begin{pro}\label{p 4.6} The decomposition
\[
W^{'}_2={\mathcal{F}}_4 \oplus {\mathcal{F}}_7
\]
is orthogonal and invariant under the action of the
group $\mathbb U(n,\mathbb R)\times Id$, \\ where
\[
\begin{array}{llll}
{\mathcal{F}}_4=\left\{F \in W^{'}_2: F(X,Y,Z)=\displaystyle{\frac{\theta (\xi)}{2n}}\{\eta(Y)
g(\varphi X,\varphi Z)-\eta(Z)g(\varphi X,\varphi Y) \}\right\}\Longleftrightarrow \\ \\
\left\{F \in W^{'}_2: {\mathcal{A}}_\xi \circ \varphi =\displaystyle{\frac{{\tr}
({\mathcal{A}}_\xi \circ \varphi )}{2n}} \, id\right\}, \\ \\
{\mathcal{F}}_7=\left\{F \in W^{'}_2: \theta (\xi)=0\right\}\Longleftrightarrow
\left\{F \in W^{'}_2: {\tr}({\mathcal{A}}_\xi \circ \varphi )=0\right\}.
\end{array}
\]
\end{pro}
Now we consider the subspace $W^{''}$ of $W_2$. We define the operators
\[
\begin{array}{llllll}
l_i: W^{''} \longrightarrow W^{''}, \quad i=1,2; \\
l_1(F)(X,Y,Z)=-\displaystyle{\frac{1}{2}}\eta(Y)\left\{g(\varphi ({\mathcal{A}} _\xi X),Z)-g(\varphi X,
{\mathcal{A}} _\xi Z)\right\}+ \\ \\
\displaystyle{\frac{1}{2}}\eta(Z)\left\{g(\varphi ({\mathcal{A}} _\xi X),Y)-g(\varphi X,
{\mathcal{A}} _\xi Y)\right\}; \\ \\
l_2(F)(X,Y,Z)=-\displaystyle{\frac{1}{2}}\eta(Y)\left\{g(\varphi ({\mathcal{A}} _\xi X),Z)+g(\varphi X,
{\mathcal{A}} _\xi Z)\right\}+ \\ \\
\displaystyle{\frac{1}{2}}\eta(Z)\left\{g(\varphi ({\mathcal{A}} _\xi X),Y)+g(\varphi X,
{\mathcal{A}} _\xi Y)\right\}.
\end{array}
\]
\begin{lem}\label{l 4.5}
The operators $l_i \, (i=1,2)$ have the following properties
\[
\begin{array}{cl}
          (i)   & l_i\circ l_i=l_i, \quad i=1,2; \\ \\

          (ii)  & \displaystyle \sum^{2}_{i=1} l_i=id; \\ \\

          (iii) & l_1\circ l_2=l_2\circ l_1=0;
\end{array}
\]
and $l_i \, (i=1,2)$ commute with $\mathbb U(n,\mathbb R)\times Id$.
\end{lem}
We denote ${\mathcal{F}}_9=Iml_1$, ${\mathcal{F}}_8=Iml_2$.
Lemma~\ref{l 4.5} implies the decomposition \\
$W^{''}={\mathcal{F}}_8 \oplus {\mathcal{F}}_9$ that is
$\mathbb U(n,\mathbb R)\times Id$ - invariant.
Using \eqref{2.7} we check that ${\mathcal{F}}_8 \bot {\mathcal{F}}_9$.
\begin{pro}\label{p 4.7} The decomposition
\[
W^{''}={\mathcal{F}}_8 \oplus {\mathcal{F}}_9
\]
is orthogonal and invariant under the action of the
group $\mathbb U(n,\mathbb R)\times Id$.
\end{pro}
Taking into account definitions of $l_1, l_2$, \eqref{4.1},
\eqref{4.3}, \eqref{2.5} we obtain
\begin{pro}\label{p 4.8}
For an arbitrary $F \in W^{''}$ we have
\[
\begin{array}{cl}
          (i)   & F=l_1F \Longleftrightarrow g({\mathcal{A}} _\xi \, . \, , \, . )=
          g( . \, , \, {\mathcal{A}} _\xi \, . ); \\
          (ii)  & F=l_1F \Longleftrightarrow F(X,Y,Z)=-F(Y,Z,X)-F(Z,X,Y); \\
          (iii) & F=l_2F \Longleftrightarrow
          g({\mathcal{A}} _\xi \, . \, , \, . )=-
          g( . \, , \, {\mathcal{A}} _\xi \, . ); \\
          (iiii)& F=l_2F \Longleftrightarrow F(X,Y,Z)=-F(Y,Z,X)+F(Z,X,Y)+2F(\varphi X,\varphi Y,Z).
\end{array}
\]
\end{pro}
From Proposition~\ref{p 4.8} the characteristic conditions of
${\mathcal{F}}_8, {\mathcal{F}}_9$ are
\begin{equation}\label{4.6}
\begin{array}{ll}
{\mathcal{F}}_8=\left\{F \in W^{''}: g({\mathcal{A}} _\xi \, . \, , \, . )=-
g( . \, , \, {\mathcal{A}} _\xi \, . )\right\} \Longleftrightarrow \\ \\
\left\{F \in W^{''}: F(X,Y,Z)=-F(Y,Z,X)+F(Z,X,Y)+2F(\varphi X,\varphi Y,Z)\right\};
\end{array}
\end{equation}
\begin{equation}\label{4.7}
\begin{array}{ll}
{\mathcal{F}}_9=\left\{F \in W^{''}: g({\mathcal{A}} _\xi \, . \, , \, . )=
g( . \, , \, {\mathcal{A}} _\xi \, . )\right\} \Longleftrightarrow \\ \\
\left\{F \in W^{''}: F(X,Y,Z)=-F(Y,Z,X)-F(Z,X,Y)\right\}.
\end{array}
\end{equation}

Finally, we denote ${\mathcal{F}}_{10}=W_3$ and
${\mathcal{F}}_{11}=W_4$. Taking into account Proposition~\ref{p
2.1}, Proposition~\ref{p 3.5}, Proposition~\ref{p 4.1},
Proposition~\ref{p 4.3}, Proposition~\ref{p 4.5},
Proposition~\ref{p 4.6}, Proposition~\ref{p 4.7} we obtain
\begin{thm}\label{4.1}
The decomposition
\[
{\mathcal{F}}={\mathcal{F}}_1 \oplus \ldots \oplus {\mathcal{F}}_{11}
\]
is orthogonal and invariant under the action of the
group $\mathbb U(n,\mathbb R)\times Id$.
\end{thm}
Next we summarize the characterization conditions for the factors ${\mathcal{F}}_i
(i=1,\ldots,11)$.
\par
Let $X,Y,Z \in V$. Then
\[
\begin{array}{llll}
{\mathcal{F}}_1: F(X,Y,Z)=\displaystyle{\frac{1}{2n}}\{g(X,\varphi Y)\theta(\varphi
Z)-g(X,\varphi Z)\theta(\varphi Y)- \\ \\
\left. g(\varphi X,\varphi Y)\theta(hZ)+g(\varphi X,\varphi Z)\theta(hY)\right\}, \\ \\
{\mathcal{F}}_2: F(\varphi X,\varphi Y,Z)=-F(X,Y,Z); \quad \theta =0, \\ \\
{\mathcal{F}}_3: F(\varphi X,\varphi Y,Z)=F(X,Y,Z), \\ \\
{\mathcal{F}}_4: F(X,Y,Z)=\displaystyle{\frac{\theta(\xi)}{2n}}\{\eta(Y)g(\varphi X,\varphi Z)-
\eta(Z)g(\varphi X,\varphi Y) \}, \\ \\
{\mathcal{F}}_5: F(X,Y,Z)=-\displaystyle{\frac{\theta^*(\xi)}{2n}}\{\eta(Y)g(X,\varphi Z)-
\eta(Z)g(X,\varphi Y) \}, \\ \\
{\mathcal{F}}_6: F(X,Y,Z)=-F(\varphi X,\varphi Y,Z)-F(\varphi X,Y,\varphi Z)=\\ \\
-F(Y,Z,X)+F(Z,X,Y)-2F(\varphi X,\varphi Y,Z); \quad \theta^*(\xi)=0, \\ \\
{\mathcal{F}}_7: F(X,Y,Z)=-F(\varphi X,\varphi Y,Z)-F(\varphi X,Y,\varphi Z)=\\ \\
-F(Y,Z,X)-F(Z,X,Y); \quad \theta(\xi)=0, \\ \\
\end{array}
\]
\[
\begin{array}{llll}
{\mathcal{F}}_8: F(X,Y,Z)=F(\varphi X,\varphi Y,Z)+F(\varphi X,Y,\varphi Z)=\\ \\
-F(Y,Z,X)+F(Z,X,Y)+2F(\varphi X,\varphi Y,Z), \\ \\
{\mathcal{F}}_9: F(X,Y,Z)=F(\varphi X,\varphi Y,Z)+F(\varphi X,Y,\varphi Z)=\\ \\
-F(Y,Z,X)-F(Z,X,Y), \\ \\
{\mathcal{F}}_{10}: F(X,Y,Z)=\eta(X)F(\xi,\varphi Y,\varphi Z), \\ \\
{\mathcal{F}}_{11}: F(X,Y,Z)=\eta(X)\{\eta(Y)\omega (Z)-\eta(Z)\omega (Y)\}.
\end{array}
\]
\section{Basic classes of almost paracontact manifolds and some examples}\label{sec-5}
\setcounter{equation}{0}
\par
Let $(M^{2n+1},\varphi,\xi,\eta,g)$ be an almost paracontact
manifold. The tensor $F$, defined by \eqref{4} we can write in the
form \eqref{2.6}, where the linear operators ${\mathcal{A}} _{e_i}
\, (i,\ldots,2n)$ and ${\mathcal{A}} _\xi$ are defined by
\[
{\mathcal{A}} _{e_i}X=(\nabla _X\varphi )e_i, \quad (i,\ldots,2n); \quad
{\mathcal{A}} _\xi X=\nabla _X\xi .
\]
We verify immediately that so defined operators ${\mathcal{A}}
_{e_i} \, (i,\ldots,2n)$ and ${\mathcal{A}} _\xi$ have the
properties \eqref{2.2}$\div$\eqref{2.5}. Using the decomposition
of the space ${\mathcal{F}}$ over $V=T_pM, \, p \in M$, we define
the corresponding subclasses of the class of almost paracontact
manifolds with respect to the covariant derivative of the
structure tensor field $\varphi$.
\par
An almost paracontact manifold is said to be in the class
${\mathcal{F}}_i \, (i=1,\ldots,11)$ if the tensor $F(X,Y,Z)=g(\nabla _X\varphi )Y,Z)$
belongs to the class ${\mathcal{F}}_i$ over $V=T_pM$ for each $p \in M$.
\par
In a similar way we define the classes ${\mathcal{F}}_i \oplus {\mathcal{F}}_j$.
It is clear that $2^{11}$ classes of almost paracontact manifolds are possible.
\par
The class ${\mathcal{F}}_0$ of almost paracontact manifolds is
defined by the condition \\ $F(X,Y,Z)=0$. This special class
belongs to everyone of the defined classes.
\begin{example}\label{5.1}
Let $(M^5,\varphi,\xi,\eta,g)$ be an almost paracontact metric
manifold. We consider a $\varphi$-basis $\{e_1,e_2,\varphi
e_1,\varphi e_2,\xi \}$ of $T_pM, \, p \in M$ such that
\[
g(e_i,e_i)=-g(\varphi e_i,\varphi e_i)=1, \quad i=1,2.
\]
\end{example}
We denote the matrixes of the the operators ${\mathcal{A}}_{e_i}$
and ${\mathcal{A}}_{\varphi e_i}$ $(i=1,2)$ with respect to the
basis $\{e_1,e_2,\varphi e_1,\varphi e_2,\xi \}$ by
${\mathcal{A}}_i$ $(i=1,2)$ and ${\mathcal{A}}_j$ $(j=3,4)$
respectively. We define
\[
{\mathcal{A}}_1=\left(
\begin{array}{rrrrr}
0 & 0 & 0 & 0 & 0 \cr a & b & c & d & 0 \cr 0 & 0 & 0 & 0 & 0 \cr
-c & -d & -a & -b & 0 \cr 0 & 0 & 0 & 0 & 0 \cr
\end{array}
\right), \quad {\mathcal{A}}_2=\left(
\begin{array}{rrrrr}
-a & -b & -c & -d & 0 \cr 0 & 0 & 0 & 0 & 0 \cr
c & d & a & b & 0
\cr 0 & 0 & 0 & 0 & 0 \cr 0 & 0 & 0 & 0 & 0 \cr
\end{array}
\right),
\]
\[
{\mathcal{A}}_3=\left(
\begin{array}{rrrrr}
0 & 0 & 0 & 0 & 0 \cr c & d & a & b & 0 \cr
0 & 0 & 0 & 0 & 0 \cr
-a & -b & -c & -d & 0 \cr 0 & 0 & 0 & 0 & 0 \cr
\end{array}
\right), \quad {\mathcal{A}}_4=\left(
\begin{array}{rrrrr}
-c & -d & -a & -b & 0 \cr
0 & 0 & 0 & 0 & 0 \cr a & b & c & d & 0 \cr
0 & 0 & 0 & 0 & 0 \cr 0 & 0 & 0 & 0 & 0 \cr
\end{array}
\right),
\]
where $a, b, c, d$ are functions over $M$.
\par
From the definitions of the matrixes ${\mathcal{A}_j}$
$(j=1,2,3,4)$ we have $\eta({\mathcal{A}}_{e_i}X)=$
\\ $\eta({\mathcal{A}}_{\varphi e_i}X)=0$ $(i=1,2)$. From
\eqref{2.4} it follows ${\mathcal{A}}_\xi X=0$. Using \eqref{2.6}
we compute
\[
F(X,Y,Z)=\left(aX^1+bX^2+cX^3+dX^4\right)\left(
Y^1Z^2-Y^2Z^1+Y^3Z^4-Y^4Z^3\right)+
\]
\[
\left(cX^1+dX^2+aX^3+bX^4\right)\left(Y^1Z^4-Y^2Z^3+Y^3Z^2-Y^4Z^1\right),
\]
where $X=X^ie_i+X^{i+2} \varphi e_i+\eta(X)\xi, \quad
Y=Y^ie_i+Y^{i+2} \varphi e_i+\eta(Y)\xi$, \\
$Z=Z^ie_i+Z^{i+2} \varphi e_i+\eta(Z)\xi \quad i=1,2$.
We verify that
\[
F(X,Y,Z)=F(\varphi X,\varphi Y,Z),
\]
which is the characterization condition of the class ${\mathcal{F}}_3$.
\begin{example}\label{5.2}
Let $(M^5,\varphi,\xi,\eta,g)$ be an almost paracontact metric manifold.
We consider a $\varphi$-basis $\{e_1,e_2,\varphi e_1,\varphi e_2,\xi \}$ of
$T_pM, \, p \in M$ such that
\[
g(e_i,e_i)=-g(\varphi e_i,\varphi e_i)=1, \quad i=1,2.
\]
\end{example}
We denote the matrixes of the operators ${\mathcal{A}}_{e_i}$,
${\mathcal{A}}_{\varphi e_i}$ $(i=1,2)$ and ${\mathcal{A}}_\xi$
with respect to the basis $\{e_1,e_2,\varphi e_1,\varphi e_2,\xi
\}$ by ${\mathcal{A}}_i$ $(i=1,2)$, ${\mathcal{A}}_j$ $(j=3,4)$
and ${\mathcal{A}}$
respectively. We define
\[
{\mathcal{A}}_1=\left(
\begin{array}{rrrrr}
0 & 0 & 0 & 0 & 0 \cr
0 & 0 & 0 & 0 & 0 \cr
0 & 0 & 0 & 0 & 0 \cr
0 & 0 & 0 & 0 & 0 \cr
-d & -e &-a &-b & 0 \cr
\end{array}
\right), \quad
{\mathcal{A}}_2=\left(
\begin{array}{rrrrr}
0 & 0 & 0 & 0 & 0 \cr
0 & 0 & 0 & 0 & 0 \cr
0 & 0 & 0 & 0 & 0 \cr
0 & 0 & 0 & 0 & 0 \cr
-e & -f & -b & -c & 0 \cr
\end{array}
\right),
\]
\[
{\mathcal{A}}_3=\left(
\begin{array}{rrrrr}
0 & 0 & 0 & 0 & 0 \cr
0 & 0 & 0 & 0 & 0 \cr
0 & 0 & 0 & 0 & 0 \cr
0 & 0 & 0 & 0 & 0 \cr
-a & -b & -d & -e & 0 \cr
\end{array}
\right), \quad
{\mathcal{A}}_4=\left(
\begin{array}{rrrrr}
0 & 0 & 0 & 0 & 0 \cr
0 & 0 & 0 & 0 & 0 \cr
0 & 0 & 0 & 0 & 0 \cr
0 & 0 & 0 & 0 & 0 \cr
-b & -c  & -e & -f & 0 \cr
\end{array}
\right),
\]
\[
{\mathcal{A}}=\left(
\begin{array}{rrrrr}
a & b & d & e & 0 \cr
b & c & e & f & 0 \cr
-d & -e & -a & -b & 0 \cr
-e & -f & -b & -c & 0 \cr
0 & 0 & 0 & 0 & 0 \cr
\end{array}
\right),
\]
where $a,b,c,d,e,f$ are functions over $M$.
Using \eqref{2.6} we compute
\[
F(X,Y,Z)=\eta(Y)\left\{a(X^1Z^3+X^3Z^1)+b(X^1Z^4+X^2Z^3+X^3Z^2+X^4Z^1)+\right.
\]
\[
c(X^2Z^4+X^4Z^2)+d(X^1Z^1+X^3Z^3)+e(X^1Z^2+X^2Z^1+X^3Z^4+X^4Z^3)+
\]
\[
\left.f(X^2Z^2+X^4Z^4)\right\}-\eta(Z)\left\{a(X^1Y^3+X^3Y^1)+
\right.
\]
\[
b(X^1Y^4+X^2Y^3+X^3Y^2+X^4Y^1)+c(X^2Y^4+X^4Y^2)+d(X^1Y^1+X^3Y^3)+
\]
\[
\left.e(X^1Y^2+X^2Y^1+X^3Y^4+X^4Y^3)+f(X^2Z^2+X^4Z^4)\right\},
\]
where $X=X^ie_i+X^{i+2} \varphi e_i+\eta(X)\xi, \quad
Y=Y^ie_i+Y^{i+2} \varphi e_i+\eta(Y)\xi$, \\
$Z=Z^ie_i+Z^{i+2} \varphi e_i+\eta(Z)\xi \quad i=1,2$.
We verify that
\[
F(X,Y,Z)=F(\varphi X,\varphi Y,Z)+F(\varphi X,Y,\varphi Z)=
-F(Y,Z,X)-F(Z,X,Y),
\]
which is the characterization condition of the class ${\mathcal{F}}_9$.

\begin{example}\label{5.3}
Let $(M^5,\varphi,\xi,\eta,g)$ be an almost paracontact metric
manifold. We consider a $\varphi$-basis $\{e_1,e_2,\varphi
e_1,\varphi e_2,\xi \}$ of $T_pM, \, p \in M$ such that
\[
g(e_i,e_i)=-g(\varphi e_i,\varphi e_i)=1, \quad i=1,2.
\]
\end{example}
We denote the matrixes of the the operators ${\mathcal{A}}_{e_i}$
and ${\mathcal{A}}_{\varphi e_i}$ $(i=1,2)$ with respect to the
basis $\{e_1,e_2,\varphi e_1,\varphi e_2,\xi \}$ by
${\mathcal{A}}_i$ $(i=1,2)$ and ${\mathcal{A}}_j$ $(j=3,4)$
respectively. We define
\[
{\mathcal{A}}_1=\left(
\begin{array}{rrrrr}
0 & 0 & 0 & 0 & 0 \cr 0 & 0 & 0 & 0 & a \cr 0 & 0 & 0 & 0 & 0 \cr
0 & 0 & 0 & 0 & b \cr 0 & 0 & 0 & 0 & 0 \cr
\end{array}
\right), \quad {\mathcal{A}}_2=\left(
\begin{array}{rrrrr}
0 & 0 & 0 & 0 & -a \cr 0 & 0 & 0 & 0 & 0 \cr 0 & 0 & 0 & 0 & -b
\cr 0 & 0 & 0 & 0 & 0 \cr 0 & 0 & 0 & 0 & 0 \cr
\end{array}
\right),
\]
\[
{\mathcal{A}}_3=\left(
\begin{array}{rrrrr}
0 & 0 & 0 & 0 & 0 \cr 0 & 0 & 0 & 0 & -b \cr 0 & 0 & 0 & 0 & 0 \cr
0 & 0 & 0 & 0 & -a \cr 0 & 0 & 0 & 0 & 0 \cr
\end{array}
\right), \quad {\mathcal{A}}_4=\left(
\begin{array}{rrrrr}
0 & 0 & 0 & 0 & b \cr 0 & 0 & 0 & 0 & 0 \cr 0 & 0 & 0 & 0 & a \cr
0 & 0 & 0 & 0 & 0 \cr 0 & 0 & 0 & 0 & 0 \cr
\end{array}
\right),
\]
where $a, b$ are functions over $M$.
\par
From the definitions of the matrixes ${\mathcal{A}_j}$ $(j=1,2,3,4)$
we have $\eta({\mathcal{A}}_{e_i}X)=$
\\ $\eta({\mathcal{A}}_{\varphi e_i}X)=0$ $(i=1,2)$. From
\eqref{2.4} it follows ${\mathcal{A}}_\xi X=0$. Using \eqref{2.6}
we compute
\[
F(X,Y,Z)=\eta(X)\left\{Y^1(aZ^2-bZ^4)+Y^2(bZ^3-aZ^1)+ \right.
\]
\[
\left.Y^3(aZ^4-bZ^2)+Y^4(bZ^1-aZ^3)\right\},
\]
where $X=X^ie_i+X^{i+2} \varphi e_i+\eta(X)\xi, \quad
Y=Y^ie_i+Y^{i+2} \varphi e_i+\eta(Y)\xi$, \\
$Z=Z^ie_i+Z^{i+2} \varphi e_i+\eta(Z)\xi \quad i=1,2$. We verify
that
\[
F(X,Y,Z)=\eta(X)F(\xi,\varphi Y,\varphi Z),
\]
which is the characterization condition of the class
${\mathcal{F}}_{10}$.
\begin{rem}\label{4.1}
Taking into account the characterization of the classes
${\mathcal{F}}_i \, (i=1,\ldots,11)$ by the linear operators
${\mathcal{A}} _{e_i} \, (i,\ldots,2n)$ and ${\mathcal{A}} _\xi$,
we can construct examples for the rest of the classes too. Using
the matrixes of the operators with respect to a $\varphi $-basis
we obtain also that an almost paracontact manifold of dimension
$3$ can not belong to the classes ${\mathcal{F}}_1,
{\mathcal{F}}_2, {\mathcal{F}}_3, {\mathcal{F}}_6$.
\end{rem}

\textbf{Acknowledgement} Simeon Zamkovoy acknowledges support from
the European Operational programm HRD through contract
BGO051PO001/07/3.3-02/53 with the Bulgarian Ministry of Education.
He also was partially supported by Contract 154/2008 with the
University of Sofia "St. Kl. Ohridski".

Galia Nakova was partially supported by Scientific researches fund
of " St. Cyril and St. Methodius " University of Veliko Tarnovo
under contract RD-491-08 from 27.06.2008.

\bibliographystyle{hamsplain}

\begin{thebibliography}{12}

\bibitem[B]{B} C. Bejan, {\em A classification of the almost parahermitian manifolds},
Proc. Conference on Diff. Geom. and Appl., Dubrovnik 1988, 23-27.

\bibitem[C]{C} V. Cruceanu, {\em Une structure parak\"{a}hlerienne sur
le fibr\'{e} tangent}, Tensor (NS) 39(1982), 81-84.

\bibitem[GB]{GB} G. Ganchev, A. Borisov, {\em Isotropic sections and curvature properties of hyperbolic K\"{a}ler manifolds},
Publ. Inst. Math. (Beograd) (N.S) 38(1985), 183-192.

\bibitem[GMG]{GMG} G. Ganchev, V. Mihova, K. Gribachev, {\em Almost Contact Manifolds with B-Metric},
Mathematica Balkanica 7(1993), 261-277.

\bibitem[GH]{GH} A. Gray, L.M. Hevella, {\em The sixteen classes of almost Hermitian manifolds and their linear invariants},
Ann. di Mat. 123(1980), 35-58.

\bibitem[MS]{MS} S. Kaneyuki, F. L. Willams, {\em Almost paracontact and parahodge structures on
manifolds}, Nagoya Math. J. 99(1985), 173-187.

\bibitem[KW]{KW} M. Manev, M. Staikova, {\em On almost paracontact Riemannian manifolds of type (n,n)}, J. Geom., 72 (2001), 108–114.

\bibitem[N]{N} A.M. Naveira, {\em A classification of Riemannian almost-product manifolds},
Rend. Mat. Appl. (7) 3(1983), 577-592.

\bibitem[Zam]{Zam} S. Zamkovoy, {\em Canonical connections on para-contact manifolds}, arXiv:0707.1787.

\end{thebibliography}






\end{document}